\documentclass[12pt,reqno]{amsart}
\usepackage{amssymb}
\usepackage{amsxtra}
\usepackage[mathscr]{eucal}

\newcommand{\Z}{{\mathbf Z}}

\newcommand{\C}{{\mathbf C}}
\newcommand{\eop}{\hfill$\square$}

\theoremstyle{plain}
\newtheorem{Thm}{Theorem}

\newtheorem{Cor}{Corollary}

\theoremstyle{definition}

\theoremstyle{remark}

\begin{document}

\title[David M.~Bradley]{On the Sum Formula for Multiple $q$-Zeta Values}

\date{\today}

\author{David~M. Bradley}
\address{Department of Mathematics \& Statistics\\
         University of Maine\\
         5752 Neville Hall
         Orono, Maine 04469-5752\\
         U.S.A.}
\email[]{bradley@math.umaine.edu, dbradley@member.ams.org}



\begin{abstract} Multiple $q$-zeta values are a 1-parameter
generalization (in fact, a $q$-analog) of the multiple harmonic
sums commonly referred to as multiple zeta values.  These latter
are obtained from the multiple $q$-zeta values in the limit as
$q\to 1$.  Here, we discuss the sum formula for multiple $q$-zeta
values, and provide a self-contained proof.  As a consequence, we
also derive a $q$-analog of Euler's evaluation of the double zeta
function $\zeta(m,1)$.
\end{abstract}

\maketitle

\interdisplaylinepenalty=500

\section{Introduction}\label{sect:Intro}

Sums of the form
\begin{equation}\label{MzvDef}
   \zeta(n_1,n_2,\dots,n_r) := \sum_{k_1>k_2>\cdots>k_r>0}\;
   \prod_{j=1}^r \frac{1}{k_j^{n_j}}
\end{equation}
have attracted increasing attention in recent years; see
eg.~\cite{BBB,BBBLa,BBBLc,BowBrad1,BowBrad3,BowBradRyoo,Prtn,DBqKarl,BK1,LeM}.
The survey articles~\cite{BowBradSurvey,Cartier,Wald,Wald2,Zud}
provide an extensive list of references. Here and throughout,
$n_1,\dots,n_r$ are positive integers with $n_1>1$, and we sum
over all positive integers $k_1,\dots,k_r$ satisfying the
indicated inequalities. Note that with positive integer arguments,
$n_1>1$ is necessary and sufficient for convergence.  The
sums~\eqref{MzvDef} are sometimes referred to as Euler sums,
because they were first studied by Euler~\cite{LE} in the case
$r=2$. In general, they may be profitably viewed as instances of
the multiple
polylogarithm~\cite{BBBLa,BowBradSurvey,Gonch1,Gonch3}, and are
now more commonly referred to as multiple zeta values, reducing to
the Riemann zeta function in the case $r=1$.  A $q$-analog
of~\eqref{MzvDef} was independently introduced
in~\cite{DBqMzv,OkudaYoshihiro,Zhao} as
\begin{equation}
   \zeta[n_1,n_2,\dots,n_r] := \sum_{k_1>k_2>\cdots >k_r>0}
   \; \prod_{j=1}^r \frac{q^{(n_j-1)k_j}}{[k_j]_q^{n_j}},
   \label{qMzvDef}
\end{equation}
where
\[
   [k]_q := \sum_{j=0}^{k-1} q^j = \frac{1-q^k}{1-q},
   \qquad 0<q<1.
\]
Observe that we now have
\[
   \zeta(n_1,\dots,n_r) = \lim_{q\to 1} \zeta[n_1,\dots,n_r],
\]
so that~\eqref{qMzvDef} represents a generalization
of~\eqref{MzvDef}.  In this note, we prove an identity
for~\eqref{qMzvDef}, the $q=1$ case of which was originally
conjectured by Moen~\cite{Hoff92} and Markett~\cite{Markett}.

It is convenient to state results in terms of the shifted multiple
zeta functions defined by
\begin{align*}
   \zeta^*[n_1,\dots,n_r] &:= \zeta[1+n_1,n_2,\dots,n_r]
   = \sum_{k_1>\cdots > k_r>0}\; \frac{q^{k_1}}{[k_1]_q}
   \prod_{j=1}^r \frac{q^{(n_j-1)k_j}}{[k_j]_q^{n_j}}\\
\intertext{and correspondingly,}
   \zeta^*(n_1,\dots,n_r) &:= \zeta(1+n_1,n_2,\dots,n_r)
   = \lim_{q\to1} \zeta^*[n_1,\dots,n_r].
\end{align*}
The main focus of our discussion is the following result.
\begin{Thm}[$q$-sum formula]\label{thm:qsum} If $N$ and $r$ are positive integers with $N\ge r$,
then
\[
   \sum_{\substack{n_1+\cdots+n_r=N\\ \forall j,\, n_j \ge 1}}
   \zeta^*[n_1,n_2,\dots,n_r] = \zeta^*[N],
\]
where the sum is over all positive integers $n_1,n_2,\dots, n_r$
such that $\sum_{j=1}^r n_j=N$.
\end{Thm}
The limiting case $q=1$ is of course the now familiar
\begin{Cor}[sum formula]\label{cor:sum} If $N$ and $r$ are positive integers with $N\ge r$,
then
\[
   \sum_{\substack{n_1+\cdots+n_r=N\\ \forall j,\, n_j \ge 1}}
   \zeta^*(n_1,n_2,\dots,n_r) = \zeta^*(N),
\]
where the sum is over all positive integers $n_1,n_2,\dots, n_r$
such that $\sum_{j=1}^r n_j=N$.
\end{Cor}
Corollary~\ref{cor:sum} was proved for $r=2$ by Euler, for $r=3$
by Hoffman and Moen~\cite{HoffMoen}, and in full generality by
Granville~\cite{Granville}.  Then Ohno derived
Corollary~\ref{cor:sum} as a consequence of his generalized
duality relation~\cite{Ohno}, and later as a consequence of an
auxiliary result used in his proof of the cyclic sum
formula~\cite{HoffOhno}.  Corollary~\ref{cor:sum} is also derived
in~\cite{OhnoZag} by specializing the height relation given there.
Subsequently and independently~\cite{DBqMzv,OkudaYoshihiro},
$q$-analogs of all these results were discovered and proved.  For
example, Theorem~\ref{thm:qsum} is derived in~\cite{DBqMzv} as a
consequence of generalized
$q$-duality~\cite{DBqMzv,OkudaYoshihiro} (a $q$-analog of the main
result in~\cite{Ohno}, but proved using an entirely different
technique).  Likewise, a $q$-analog of the cyclic sum
formula~\cite{DBqMzv,OkudaYoshihiro} also leads to a quick
proof~\cite{DBqMzv} of Theorem~\ref{thm:qsum}. Finally,
in~\cite{OkudaYoshihiro}, a $q$-analog of the height relation is
also given; we show below how this too can be used to derive
Theorem~\ref{thm:qsum}.  However, as all these proofs of
Theorem~\ref{thm:qsum} depend on comparatively more sophisticated
results for~\eqref{qMzvDef}, we feel it may be of interest to give
a self-contained proof, more in the spirit of~\cite{Granville}.

\section{Self-Contained Proof of Theorem~\ref{thm:qsum}}\label{sect:self}

By expanding both sides in powers of $z$ and comparing
coefficients, one readily sees that Theorem~\ref{thm:qsum} is
equivalent to the following result.

\begin{Thm}\label{thm:qsumgf} If $r$ is a positive integer and $z\in\C\setminus\{q^{-m}[m]_q :
m\in\Z^{+}\}$, then
\begin{equation}\label{qsumgf}
   \sum_{k_1>\cdots >k_r>0}\; \frac{q^{k_1}}{[k_1]_q}
   \prod_{j=1}^r \frac{1}{[k_j]_q-zq^{k_j}}
   = \sum_{m=1}^\infty
   \frac{q^{rm}}{[m]_q^r\left([m]_q-zq^m\right)}.
\end{equation}
\end{Thm}

\noindent{\bf Proof of Theorem~\ref{thm:qsumgf}.} Let $L_r
=L_r(z)$ denote the left hand side of~\eqref{qsumgf}.  By partial
fractions,
\begin{equation}\label{Lr=sumSj}
   L_r =\sum_{j=1}^r S_j 
\end{equation}
where
\[
   S_j = S_{j,r}(z)
       := \sum_{k_1>\cdots >k_r>0}\;\frac{q^{k_1}}{[k_1]_q\left([k_j]_q-zq^{k_j}\right)}
          \prod_{\substack{i=1\\i\ne j}}^r\frac{1}{[k_i-k_j]_q}.
\]
Now rename $k_j=m$ and sum first on $m$, so that
\begin{equation}\label{Sj}
   S_j = \sum_{m=1}^\infty \frac{A(m,j-1)B(m,r-j)}{[m]_q-zq^m} ,
\end{equation}
where $A(m,0) := q^{m}/[m]_q$,
\[
   A(m,j-1) := \sum_{k_1>\cdots >k_{j-1}>m}\;
   \frac{q^{k_1}}{[k_1]_q} \prod_{i=1}^{j-1}\frac{1}{[k_i-m]_q}
   \quad\mbox{for}\quad 2\le j\le r,
\]
$B(m,0) := 1$ and for $1\le j\le r-1$,
\[
   B(m,r-j)
   := \sum_{m>k_{j+1}>\cdots > k_r>0}\; \prod_{i=j+1}^r \frac{1}{[k_i-m]_q}
    = (-1)^{r-j}\sum_{m>k_{j+1}>\cdots > k_r>0}\; \prod_{i=j+1}^r \frac{q^{m-k_i}}{[m-k_i]_q}.
\]
From~\eqref{Lr=sumSj} and~\eqref{Sj} we now get that
\[
   L_r
   =\sum_{j=0}^{r-1} S_{j+1}
   = \sum_{m=1}^\infty\frac{1}{[m]_q-zq^m}\sum_{j=0}^{r-1} A(m,j)B(m,r-1-j),
\]
and hence
\begin{equation}\label{Lrgf}
   \sum_{r=1}^\infty x^{r-1} L_r = \sum_{m=1}^\infty
   \frac{A_m(x)B_m(x)}{[m]_q-zq^m},
\end{equation}
where the generating functions $A_m$ and $B_m$ are defined by
\[
   A_m(x) := \sum_{n=0}^\infty x^n A(m,n),\qquad
   B_m(x) := \sum_{n=0}^\infty x^n B(m,n).
\]
The proof of Theorem~\ref{thm:qsumgf} now follows more or less
immediately from the representations
\begin{equation}\label{ABreps}
   A_m(x) = \frac{q^m}{[m]_q}\prod_{c=1}^m
   \bigg(1-\frac{xq^c}{[c]_q}\bigg)^{-1} \quad\mbox{and}\quad
   B_m(x) = \prod_{b=1}^{m-1} \bigg(1-\frac{xq^b}{[b]_q}\bigg).
\end{equation}
To see this, observe that~\eqref{ABreps} gives
\[
   A_m(x)B_m(x)
   =\frac{q^m}{[m]_q}\bigg(1-\frac{xq^m}{[m]_q}\bigg)^{-1}
   =\sum_{r=1}^\infty x^{r-1}\frac{q^{rm}}{[m]_q^r},
\]
and hence from~\eqref{Lrgf},
\[
   \sum_{r=1}^\infty x^{r-1}L_r
   = \sum_{r=1}^\infty x^{r-1}\sum_{m=1}^\infty\frac{q^{rm}}{[m]_q^r\left([m]_q-zq^m\right)}.
\]
It now remains only to prove the representations~\eqref{ABreps}.
First, note that
\begin{align*}
   B_m(x)
   &= \sum_{n=0}^\infty x^n (-1)^n \sum_{m>k_1>\cdots >k_n>0}\; \prod_{j=1}^n \frac{q^{m-k_j}}{[m-k_j]_q}
    = \sum_{n=0}^\infty (-x)^n \sum_{m>b_n>\cdots >b_1>0}\;\prod_{j=1}^n \frac{q^{b_j}}{[b_j]_q}\\
   &= \prod_{b=1}^{m-1}\bigg(1-\frac{xq^b}{[b]_q}\bigg).
\end{align*}
Next, we define
\[
   A(m,n,k) := \sum_{b_1>\cdots >
   b_n>k}\;\frac{q^{m+b_1}}{[m+b_1]_q}\prod_{j=1}^n\frac{1}{[b_j]_q},
\]
and note that $A(m,n) = A(m,n,0)$.  We have
\[
   A(m,1,k) = \sum_{b>k} \frac{q^{m+b}}{[m+b]_q[b]_q}
   =
   \frac{q^m}{[m]_q}\sum_{b>k}\bigg(\frac{q^b}{[b]_q}-\frac{q^{m+b}}{[m+b]_q}\bigg)
   =\frac{q^m}{[m]_q}\sum_{m\ge c\ge 1} \frac{q^{c+k}}{[c+k]_q},
\]
and if for some positive integer $n$,
\[
   A(m,n,k)
   = \frac{q^m}{[m]_q}\sum_{m\ge c_1\ge\cdots\ge c_n\ge1}\;\frac{q^{c_n+k}}{[c_n+k]_q}\prod_{j=1}^{n-1}\frac{q^{c_j}}{[c_j]_q},
\]
then
\begin{align*}
   A(m,n+1,k)
   &= \sum_{b_1>\cdots > b_{n+1}>k}\;\frac{q^{m+b_1}}{[m+b_1]_q}\prod_{j=1}^{n+1}\frac{1}{[b_j]_q}\\
   &= \sum_{b_2>\cdots >b_{n+1}>k}\;\bigg(\prod_{j=2}^{n+1}\frac{1}{[b_j]_q}\bigg)\sum_{b_1>b_2}\frac{q^{m+b_1}}{[m+b_1]_q[b_1]_q}\\
   &= \sum_{b_2>\cdots >b_{n+1}>k}\;\bigg(\prod_{j=2}^{n+1}\frac{1}{[b_j]_q}\bigg)A(m,1,b_2)\\
   &= \sum_{b_2>\cdots >b_{n+1}>k}\;\bigg(\prod_{j=2}^{n+1}\frac{1}{[b_j]_q}\bigg)
      \frac{q^m}{[m]_q}\sum_{c_0=1}^m\frac{q^{c_0+b_2}}{[c_0+b_2]_q}\\
   &= \frac{q^m}{[m]_q}\sum_{c_0=1}^m\;\sum_{b_2>\cdots>b_{n+1}>k}\;\frac{q^{c_0+b_2}}{[c_0+b_2]_q}\prod_{j=2}^{n+1}\frac{1}{[b_j]_q}\\
   &= \frac{q^m}{[m]_q}\sum_{c_0=1}^m A(c_0,n,k)\\
   &= \frac{q^m}{[m]_q}\sum_{m\ge c_0\ge \cdots \ge c_n\ge 1}\;\frac{q^{c_n+k}}{[c_n+k]_q}\prod_{j=0}^{n-1}\frac{q^{c_j}}{[c_j]_q},
\end{align*}
by the induction hypothesis.   It follows that
\[
   A(m,n) = A(m,n,0) = \frac{q^m}{[m]_q}\sum_{m\ge c_1\ge\cdots\ge c_n\ge 1}\;\prod_{j=1}^n \frac{q^{c_j}}{[c_j]_q},
\]
and hence
\[
   A_m(x) =
   \frac{q^m}{[m]_q}\prod_{c=1}^m
   \bigg(1+\frac{xq^c}{[c]_q}+\bigg(\frac{xq^c}{[c]_q}\bigg)^2+\bigg(\frac{xq^c}{[c]_q}\bigg)^3+\cdots\bigg)
  = \frac{q^m}{[m]_q}\prod_{c=1}^m\bigg(1-\frac{xq^c}{[c]_q}\bigg)^{-1}.
\]
\eop

\section{Evaluation of $\zeta[m,1]$}

Euler~\cite{LE,Niels} (see also~\cite[eq.~(31)]{BBB}) proved that
for all integers $m\ge 2$,
\[
   2\zeta(m,1) = m\zeta(m+1) -
   \sum_{k=1}^{m-2}\zeta(m-k)\zeta(k+1),
\]
thereby expressing $\zeta(m,1)$ in terms of values of the Riemann
zeta function.  The following $q$-analog of Euler's formula is an
easy consequence of the $r=2$ case of Theorem~\ref{thm:qsum} and
the $q$-stuffle multiplication rule~\cite{DBqMzv}.

\begin{Cor}[Corollary 8 of~\cite{DBqMzv}]\label{cor:qEuler} Let $2\le m\in\Z$.  Then
\[
   2\zeta[m,1] = m\zeta[m+1] +(1-q)(m-2)\zeta[m] -
   \sum_{k=1}^{m-2}\zeta[m-k]\,\zeta[k+1].
\]
\end{Cor}
In particular, when $m=2$ we
get $\zeta[2,1]=\zeta[3]$, which is
probably the simplest non-trivial identity satisfied by the
multiple $q$-zeta function.

\noindent{\bf Proof.} For $1\le k\le m-2$ the $q$-stuffle
multiplication rule~\cite{DBqMzv} implies that
\[
   \zeta[m-k]\zeta[k+1]
   =\zeta[m+1]+(1-q)\zeta[m]+\zeta[m-k,k+1]+\zeta[k+1,m-k].
\]
Summing on $k$, we find that
\[
   \sum_{k=1}^{m-2} \zeta[m-k]\zeta[k+1]
   =(m-2)\left(\zeta[m+1]+(1-q)\zeta[m]\right)
      +2\sum_{\substack{s+t=m+1\\s,t\ge 2}}\zeta[s,t].
\]
But Theorem~\ref{thm:qsum} gives
\[
   \sum_{\substack{s+t=m+1\\s,\, t\ge 2}}\zeta[s,t]
   =\sum_{\substack{s+t=m+1\\s\ge 2,\, t\ge 1}}\zeta[s,t]-\zeta[m,1]
   =\zeta[m+1]-\zeta[m,1].
\]
It follows that
\[
   \sum_{k=1}^{m-2} \zeta[m-k]\zeta[k+1]
   = m\zeta[m+1]+(1-q)(m-2)\zeta[m]-2\zeta[m,1].
\]
\eop

\section{Height Relation}\label{sect:height}

Corollary~\ref{cor:qEuler} is also derived in~\cite{DBqMzv} as a
consequence of the more general double generating function
identity
\begin{multline}\label{Drin}
   \sum_{m=0}^\infty\sum_{n=0}^\infty x^{m+1}y^{n+1}
   \zeta[m+2,\{1\}^n]\\ = 1-\exp\bigg\{\sum_{k=2}^\infty
   \big\{x^k+y^k-\big(x+y+(1-q)xy\big)^k\big\}\frac1k\sum_{j=2}^k
   (q-1)^{k-j}\zeta[j]\bigg\},
\end{multline}
which implies, among other things, that
$\zeta[m+2,\{1\}^n]=\zeta[n+2,\{1\}^m]$ can be expressed in terms
of sums of products of single $q$-zeta values for every pair of
non-negative integers $m$ and $n$.  In fact~\eqref{Drin} is just
the constant term of an even more general result.

For any multi-index $\vec n=(n_1,\dots,n_r)$ of positive integers,
the weight, depth, and height of $\vec n$ are the integers
$n=n_1+n_2+\cdots+n_r$, $r$, and $s=\#\{j:n_j>1\}$, respectively.
Denote the set of multi-indices of weight $n$, depth $r$ and
height $s$ with the additional requirement $n_1>1$ by
$I_0(n,r,s)$, and set
\[
   G_0[n,r,s]:= \sum_{\vec n\in I_0(n,r,s)} \zeta[\vec n],
   \qquad
   \Phi_0[x,y,z]:= \sum_{n,r,s=0}^\infty G_0[n,r,s]
   x^{n-r-s}y^{r-s} z^{s-1}.
\]
Okuda and Takeyama~\cite{OkudaYoshihiro} proved that
\begin{align*}
   1+(z-xy)\Phi_0[x,y,z] &= \prod_{n=1}^\infty \frac{([n]_q-\alpha
   q^n)([n]_q-\beta q^n)}{([n]_q-xq^n)([n]_q-yq^n)}\\
   &= \exp\bigg\{\sum_{k=2}^\infty
   \big(x^k+y^k-\alpha^k-\beta^k\big)\frac1k\sum_{j=2}^k
   (q-1)^{k-j}\zeta[j]\bigg\},
\end{align*}
where $\alpha$ and $\beta$ are determined by
\[
   \alpha +\beta = x+y+(q-1)(z-xy),
   \qquad
   \alpha\beta = z.
\]
The limiting case $q\to 1$ reduces to the height relation
of~\cite{OhnoZag}. The case $z=0$ gives~\eqref{Drin}.  As with the
$q=1$ case~\cite{OhnoZag}, taking the limit as $z\to xy$ gives
\[
   \Phi_0[x,y,xy] = \sum_{m=1}^\infty
   \frac{q^m}{([m]_q-xq^m)([m]_q-yq^m)}
   = \sum_{n>r>0} \zeta[n] x^{n-r-1}y^{r-1}.
\]
On the other hand, by definition,
\[
   \Phi_0[x,y,xy] = \sum_{n>r>0} G_0[n,r] x^{n-r-1}y^{r-1},
\]
where $G_0[n,r]$ is the sum of all multiple $q$-zeta values of
weight $n$ and depth $r$.  Thus, we obtain $G_0[n,r]=\zeta[n]$
i.e.\ Theorem~\ref{thm:qsum} again. \qed

\end{document}